# BOUNDING MINIMUM DISTANCES OF CYCLIC CODES USING ALGEBRAIC GEOMETRY

Nigel Boston

**0. Overview.**

There are many results on the minimum distance of a cyclic code of the form that if a certain set $T$ is a subset of the defining set of the code, then the minimum distance of the code is greater than some integer $t$. This includes the BCH, Hartmann-Tzeng, Roos, and shift bounds and generalizations of these (see below). In this paper we define certain projective varieties $V(T, t)$ whose properties determine whether, if $T$ is in the defining set, the code has minimum distance exceeding $t$. Thus our attention shifts to the study of these varieties. By investigating them we will prove various new bounds. It is interesting, however, to note that there are cases that existing methods handle, that our methods do not, and vice versa. We end with a number of conjectures.

**1. Some Definitions.**

A (linear) code $C$ is a subspace of a finite-dimensional vector space $V$ over a finite field $\mathbf{F}_q$. The code has parameters $n := \dim V, k := \dim C, d := \min\{w(x) : 0 \neq x \in C\}$ where $w(x)$ is the number of nonzero components of $x$ (with respect to a fixed basis of $V$). Such a code can correct $t$ errors if $t < d/2$, i.e. given $x \in V$ (the received word) there is at most one codeword (element of $C$) that differs from $x$ in $< d/2$ places. One measure of a good code is to have $d$ large. More precisely, we want the rate $k/n$ and the relative minimum distance $d/n$ to both be large. These are competing requirements.

A source of good codes of interest to us is the family of cyclic codes. Let $V = \mathbf{F}_q[x]/(x^n - 1)$ explicitly considered as a vector space over $\mathbf{F}_q$ by the map $V \to \mathbf{F}_q^n$ sending $a_0 + a_1 x + ... + a_{n-1} x^{n-1} \mapsto (a_0, a_1, ..., a_{n-1})$. Suppose that $g \in \mathbf{F}_q[x]$ is monic and divides $x^n - 1$ and that the degree of $g$ is $n - k$. Then the ideal of $V$ generated by $g$ yields a subspace of $V$ represented by $C = \{h(x)g(x) : \deg h \leq k-1\}$. This then is a code with parameters $[n, k, ?]$.

An example of this is the (binary) Golay code. For this, we take $n = 23$ and $S$ to be the nonzero squares modulo 23. Since 23 divides $2^{11} - 1$, there exists $\alpha$ of order 23 in $\mathbf{F}_{2^{11}}^*$. Let $g(x) = \prod_{i \in S}(x - \alpha^i)$. Since $g \in \mathbf{F}_2[x]$, the construction of the previous paragraph yields a code and its parameters are $[23, 12, 7]$.

What do we know of the minimum distance of a cyclic code? First of all, there is the BCH bound. We will, for the rest of this paper, fix the notation of the code $C$, the ambient vector space $V$, the finite field $\mathbf{F}_q$, and the polynomial $g$ producing $C$. Let us also fix $\alpha \in \overline{\mathbf{F}_q}^*$ of order $n$ (so we are assuming $n$ and $q$ to be coprime). Since $g$ is monic and divides $x^n - 1$, it factors in $\overline{\mathbf{F}_q}[x]$ as $\prod_{i \in S}(x - \alpha^i)$ for some $S \subseteq \{0, 1, 2, ..., n-1\}$. This set $S$ is called the defining set of the code $C$.

**Theorem.** (Bose-Ray Chaudhuri-Hocquenghem) If $\{r, r+1, ..., r+t-1\} \subseteq S$, then the code has minimum distance $d > t$.

This theorem is usually called the BCH bound. For example, for the binary Golay code, since $\{1, 2, 3, 4\} \subseteq S$, its minimum distance $d > 4$.

We now recall a proof of the BCH bound by a method that generalizes to yield several results of the form: if $S$ contains a certain set, then this implies a lower bound on the minimum distance.

*Proof.* Suppose that $c(x) = h(x)g(x)$ is a codeword, say it's $\sum_{i=0}^{n-1} c_i x^i$. The hypothesis on $S$ implies that $c(\alpha^{r+j}) = 0$ for $0 \leq j \leq t-1$. If $c(x)$ has at most $t$ nonzero coefficients, we let $c_{i_1}, ..., c_{i_t}$ include them. Plugging in, we get a matrix equation, where the determinant of the matrix is a power of $\alpha$ times a Vandermonde determinant so is nonzero. Thus, the only solution is $c_{i_1} = ... = c_{i_t} = 0$ so that $c(x)$ is the zero codeword. In other words, no other codeword has weight $\leq t$.

We will also be interested in other known bounds, such as the following.

**Theorem (Hartmann-Tzeng).** Suppose $\{r, r+1, ..., r+t-1, r+m, r+m+1, ..., r+m+t-1, r+2m, r+2m+1, ..., r+2m+t-1, ..., r+km+t-1\} \subseteq S$, and $m$ is relatively prime to $n$, then $d > t + k$.

The shortcoming of this and similar results is the special form of the subset of $S$. In practice, $S$ will not have such a pleasant subset. Consider, for instance, quadratic residue codes, for which $S$ is the set of nonzero squares modulo some prime. The aim of this paper is to answer the question of what we can conclude in such situations.

**2. The Varieties $V(T,t)$.**

We begin by defining the projective varieties of interest to us and then relate them to cyclic codes. We end this section by listing some elementary properties of the varieties.

For any set $U$ of $m$ integers we use $\Delta[U]$ to denote the determinant of the $m$ by $m$ matrix whose rows are of the form $(x_1^i, x_2^i, ...x_m^i)$ as $i$ ranges over $U$. (The $x_1, ..., x_m$ should be considered as dummy variables.) Thus, if e.g. $U = \{0, 1, 2, ..., m-1\}$, then $\Delta[U]$ is simply the Vandermonde determinant $\Delta_m$. Furthermore, $f[U]$ will denote the homogeneous polynomial $\Delta[U]/\Delta_m$. Finally, $V(T,t)$ will denote the common zero set of the $f[U]$ as $U$ runs over all subsets $U$ of $T$ of cardinality $t$, considered as a subvariety of $\mathbf{P}^{t-1}$. That is,

$$V(T,t) := \{\underline{x} = (x_1 : x_2 : ... : x_t) \in \mathbf{P}^{t-1} : f[U](\underline{x}) = 0 \text{ for all } U \subseteq T, \#U = t\}$$

A first simple example of this (many more will be given in the next section) is the case $T = \{0, 1, 2, ..., t-1\}$. Here, $f[T] = 1$ and so $V(T, t) = \emptyset$. We can now generalize the above proof of the BCH bound.

**Theorem.** If $T$ is a subset of the defining set of the cyclic code $C$ and if $V(T, t)$ has no points of the form $(\alpha^{i_1} : \alpha^{i_2} : ... : \alpha^{i_t})$ ($i_1, ..., i_t$ distinct mod $n$), then the code $C$ has minimum distance greater than $t$.

*Proof.* Suppose that $c(x) = h(x)g(x)$ is a codeword, say it's $\sum_{i=0}^{n-1} c_i x^i$. Let $U$ be any subset of $T$ of cardinality $t$. The hypothesis on the defining set of the code implies that $c(\alpha^i) = 0$ for $i \in U$. If $c(x)$ has at most $t$ nonzero coefficients, we let $c_{i_1}, ..., c_{i_t}$ include them. Plugging in, we get a matrix equation, where the determinant of the matrix is $\Delta[U]$ evaluated at $x_j = \alpha^{i_j}$. Since the $i_j$ are distinct modulo $n$, this determinant is zero if and only if $f[U]$ is zero at $(\alpha^{i_1} : \alpha^{i_2} : ... : \alpha^{i_t})$. Our hypothesis on $V(T, t)$ implies that this cannot happen for all such $U$. Thus we get a contradiction and so each codeword has $> t$ nonzero entries.

**Corollary.** Since $V(\{r, r+1, ..., r+t-1\}, t) = \emptyset$, the BCH bound follows.

We note a few elementary results concerning these varieties. First, $V(T, t) = V(T+m, t)$ since the determinants involved differ only by powers of $\alpha$, which are never zero (yielding e.g. the above corollary). Thus, we will typically assume that the smallest element of $T$ is 0. Second, in computing $V(T, t)$ we do not need to consider all $t$-element subsets $U$ of $T$. For instance, to compute $V(\{0, 1, 3, 4\}, 3)$, we only need intersect $V(\{0, 1, 3\}, 3)$ and $V(\{0, 1, 4\}, 3)$.

## 3. Examples of $V(T, t)$.

(1) $T = \{0, 1, 3, 4\}$, $t = 3$. Then $V(T, t)$ is the intersection of curves $x_1 + x_2 + x_3 = 0$ and $x_1^2 + x_2^2 + x_3^2 + x_2 x_3 + x_3 x_1 + x_1 x_2 = 0$. Substituting for $x_3 = -x_1 - x_2$ yields $x_1^2 + x_1 x_2 + x_2^2 = 0$ and so $V(T, t) = \{(1 : \omega : \omega^2), (1 : \omega^2 : \omega)\}$, where $\omega$ is a primitive cube root of 1.

**Corollary.** If 3 does not divide $n$ and $\{0, 1, 3, 4\} \subseteq S$, then $d \geq 4$.

*Proof.* If 3 does not divide $n$, then $V(T, t)$ contains no points of the form $(\alpha^{i_1} : \alpha^{i_2} : \alpha^{i_3})$ since $\alpha_{i_1 - i_2}$ cannot be a cube root of 1.

(2) $T = \{0, 1, 3, 5\}$, $t = 3$. Then $V(T, t)$ is the intersection of curves $x_1 + x_2 + x_3 = 0$ and $V(\{0, 1, 5\}, 3) = \{x_1^3 + x_2^3 + x_3^3 + x_1^2 x_2 + x_1 x_2^2 + x_1^2 x_3 + x_1 x_3^2 + x_2^2 x_3 + x_2 x_3^2 + x_1 x_2 x_3 = 0\}$. Substituting for $x_3 = -x_1 - x_2$ yields $x_1 x_2 (x_1 + x_2) = 0$ and so $V(T, t) = \{(0 : 1 : -1), (1 : 0 : -1), (1 : -1 : 0)\}$.

**Corollary.** If $\{0, 1, 3, 5\} \subseteq S$, then $d \geq 4$.

*Proof.* $V(T, t)$ contains no points of the form $(\alpha^{i_1} : \alpha^{i_2} : \alpha^{i_3})$ since $\alpha \neq 0$.

(3) $T = \{0, 1, 3, 8\}$, $t = 3$. Then $V(T, t)$ is the intersection of curves $x_1 + x_2 + x_3 =$

0 and $V(\{0,1,8\},3)$, a curve of degree 6. Substituting for $x_3 = -x_1 - x_2$ yields $x_1^6 + 3x_1^5 x_2 + 7x_1^4 x_2^2 + 9x_1^3 x_2^3 + 7x_1^2 x_2^4 + 3x_1 x_2^5 + x_2^6 = 0$. Let $\beta$ be a root of the irreducible (over $\mathbf{Q}$) polynomial $u^6 + 3u^5 + 7u^4 + 9u^3 + 7u^2 + 3u + 1$. Then $V(T,t)$ consists of the 6 points obtained by permuting the coordinates of $(1 : \beta : -1 - \beta)$.

To analyze $V(T,t)$ further, we need to know what the field $K = \mathbf{Q}(\beta)$ is. Amazingly, it turns out to be a well-known field, namely the Hilbert class field of $\mathbf{Q}(\sqrt{-31})$. This is proven by computing the discriminant of $K$ to be $-31^3$ and noting that there is only one $S_3$-extension of $\mathbf{Q}$ ramified only at 31. (It is clear that $K$ is Galois over $\mathbf{Q}$ with Galois group $S_3$ since $V(T,t)$ is defined by equations symmetric in $x_1, x_2, x_3$.)

Note that the smallest field of characteristic $p$ over which the points of $V(T,t)$ are defined is $\mathbf{F}_p(\beta)$ (for any $p$). These fields are precisely the residue fields of $K$. These residue fields are nicely described by class field theory - for instance, in characteristic $p$ it equals $\mathbf{F}_p$ if and only if $p$ is of the form $u^2 + 31v^2$ for some integers $u, v$.

**Corollary.** If $\{0,1,3,8\} \subseteq S$ and $\mathbf{F}_q(\alpha)$ does not contain a residue field of $K$, then $d \geq 4$.

*Proof.* $V(T,t)$ contains no points of the form $(\alpha^{i_1} : \alpha^{i_2} : \alpha^{i_3})$ since $\mathbf{F}_q(\alpha)$ does not contain the field of definition of $V(T,t)$ in that characteristic.

(4) $T = \{0, 1, 3, 14\}, t = 3$. Then $V(T,t)$ is the intersection of curves $x_1 + x_2 + x_3 = 0$ and $V(\{0, 1, 14\}, 3)$, a curve of degree 12. Substituting for $x_3 = -x_1 - x_2$ yields $F(x_1, x_2) = 0$ where $F(u, 1)$ is the irreducible polynomial $u^{12} + 6u^{11} + 31u^{10} + 100u^9 + 221u^8 + 350u^7 + 407u^6 + 350u^5 + 221u^4 + 100u^3 + 31u^2 + 6u + 1$. The Galois group $G$ of this is TransitiveGroup(12,35) in the MAGMA database, of order 72 (current computer algebra systems, e.g. KASH, can compute Galois groups up to degree 15). Then $V(T,t)$ consists of the 12 points $\{(-1 : \beta : -1 - \beta) : F(u,1)(\beta) = 0\}$.

As before, let $K = \mathbf{Q}(\beta)$. Its residue fields are harder to describe than in the previous example, but it is worth noting that its Galois group $G$ is solvable and so class field theory gives some results on the explicit form of these residue fields (although unfortunately the group is not even metabelian). In fact $G$ is a semidirect product of the dihedral group of order 8 by a normal elementary abelian subgroup of order 9.

**Corollary.** If $\{0, 1, 3, 14\} \subseteq S$ and $\mathbf{F}_q(\alpha)$ does not contain a residue field of $K$, then $d \geq 4$.

*Proof.* As in the previous example.

(5) $T = \{0, 1, 3, 4, 6\}, t = 4$. Then $V(T,t)$ is the intersection of surfaces $x_1 x_2 + x_1 x_3 + x_1 x_4 + x_2 x_3 + x_2 x_4 + x_3 x_4 = 0$ and $x_1^3 x_2 + x_1^2 x_2^2 + x_1 x_2^3 + x_1^3 x_3 + 2x_1^2 x_2 x_3 + 2x_1 x_2^2 x_3 + x_2^3 x_3 + x_1^2 x_3^2 + 2x_1 x_2 x_3^2 + x_2^2 x_3^2 + x_1 x_3^3 + x_2 x_3^3 + x_1^3 x_4 + 2x_1^2 x_2 x_4 + 2x_1 x_2^2 x_4 + x_2^3 x_4 + 2x_1^2 x_3 x_4 + 3x_1 x_2 x_3 x_4 + 2x_2^2 x_3 x_4 + 2x_1 x_3^2 x_4 + x_3^3 x_4 + x_1^2 x_4^2 + 2x_2 x_3^2 x_4 + 2x_1 x_2 x_4^2 + x_2^2 x_4^2 + 2x_1 x_3 x_4^2 + 2x_2 x_3 x_4^2 + x_3^2 x_4^2 + x_1 x_4^3 + x_2 x_4^3 + x_3 x_4^3 = 0$. (Note: these equations are obtained easily using Mathematica.) This is the union of 8 lines, namely $\{x_3 = \omega x_2 = \omega^2 x_1\}$ together with the other 7 lines obtained by permuting coordinates ($\omega$ as usual a primitive cube root of 1). Explicit calculation shows that this is the same variety as $V(\{0, 1, 3, 4, 6, 7\}, 4)$. We therefore get a corollary stronger than what Hartmann-Tzeng gives us, but in fact the varieties are trivially

equal (see next section).

The variety $V(T,t)$ is not a finite set of points this time, but we can still obtain applications to coding theory.

**Corollary.** If $\{0,1,3,4,6\} \subseteq S$ and 3 does not divide $n$, then $d \geq 5$.

*Proof.* Since the ratio of some of the coordinates of a point in $V(T,t)$ must be a cube root of 1, $V(T,t)$ contains no point of the form $(\alpha^{i_1} : \alpha^{i_2} : \alpha^{i_3} : \alpha^{i_4})$.

(6) $T = \{0,1,2,4,5,6,8\}, t = 5$. Then $V(T,t)$ is the intersection of 3 hypersurfaces and one can calculate that it consists of 30 lines, where each line is given by having one of the 5 variables be arbitrary and the remaining 4 variables be $1, \zeta^{i_1}, \zeta^{i_2}, \zeta^{i_3}$, where $\zeta$ is a primitive 4th root of 1 and $0, i_1, i_2, i_3$ are all distinct modulo 4. Explicit calculation shows that this is the same variety as $V(\{0,1,2,4,5,6,8,9,10\},5)$. We therefore get a corollary stronger than what Hartmann-Tzeng gives us, but in fact the varieties are trivially equal (see next section).

The variety $V(T,t)$ is not a finite set of points this time, but we can still obtain ap plications to coding theory.

**Corollary.** If $\{0,1,2,4,5,6,8\} \subseteq S$ and 4 does not divide $n$ (note that this is a stronger result than the Hartmann-Tzeng condition that $n$ should be relatively prime to 4), then $d \geq 6$.

*Proof.* As with the previous example.

(7) $T = \{0,1,3,4,6,7\}, t = 5$. Then $V(T,t)$ is the intersection of 2 hypersurfaces and one can calculate that it consists of 20 planes, where each plane is given by having 2 of the 5 variables be arbitrary and the remaining 3 variables be $1, \omega^{i_1}, \omega^{i_2}$ where $\omega$ is a primitive cube root of 1 and $0, i_1, i_2$ are all distinct modulo 3. Explicit calculation shows that this is the same variety as $V(\{0,1,3,4,6,7,9,10\},5)$. We therefore get a corollary stronger than what Hartmann-Tzeng gives us, but in fact the varieties are trivially equal (see next section). The corollary below is, however, stronger than Hartmann-Tzeng and stronger than what more elementary coding theory methods yield.

The variety $V(T,t)$ is not a finite set of points this time, but we can still obtain applications to coding theory.

**Corollary.** If $\{0,1,3,4,6,7\} \subseteq S$ and 3 does not divide $n$, then $d \geq 6$.

*Proof.* As with the previous example.

(8) $T = \{0,1,4,5,8\}, t = 4$. Then $V(T,t)$ is the intersection of 2 surfaces and one can calculate that it consists of 24 lines of a predictable form (see the next section) together with a conic, namely the intersection of $x_1 + x_2 + x_3 + x_4 = 0$ and $x_2^2 + x_3^2 + x_4^2 + x_3x_4 + x_4x_2 + x_2x_3 = 0$. Explicit calculation shows that this is the same variety as $V(\{0,1,4,5,8,9\},4)$. Using that conics have rational parametrizations, it follows that $V(T,t)$ only has points of the desired form if 4 divides $n$.

**Corollary.** If $\{0,1,4,5,8\} \subseteq S$ and 4 does not divide $n$, then $d \geq 5$.

(9) $T = \{0,1,5,6,10\}, t = 4$. Then $V(T,t)$ is the intersection of 2 surfaces and one can calculate that it consists of 48 lines of a predictable form (see the next section) together with a curve of genus 4. This last curve is the intersection of

quadric surface $x_1^2 + x_2^2 + x_3^2 + x_4^2 + x_1x_2 + x_1x_3 + x_1x_4 + x_2x_3 + x_2x_4 + x_3x_4 = 0$ and $x_2^3 + x_3^3 + x_4^3 + x_2^2x_3 + x_2x_3^2 + x_2^2x_4 + x_2x_4^2 + x_3x_4^2 + x_3^2x_4 + x_2x_3x_4 = 0$. This latter surface is $\mathbf{P}^1 \times E$, where $E$ is the elliptic curve denoted $50A$ in Cremona's tables. In fact, letting the curve of genus 4 be denoted $X$, we get that its Jacobian is a product of 4 copies of $E$, so that if $\#E(\mathbf{F}_q) = q+1-a_q$, then $\#X(\mathbf{F}_q) = q+1-4a_q$ ($q$ not a power of 2 or 5). It is harder now to ensure that $V(T,t)$ has no points over $\mathbf{F}_q$) for a given $n$ since the points of $X(\mathbf{F}_q)$ are harder to describe. We can, however, get results such as:

**Corollary.** If $\mathbf{F}_q(\alpha) = \mathbf{F}_{47}$, then if $\{0, 1, 5, 6, 10\} \subseteq S$, then $d \geq 5$.

*Proof.* $a_{47} = 12$ and so $\#X(\mathbf{F}_q) = 0$. Moreover, there are no primitive 5th roots of 1 in $\mathbf{F}_{47}$ and so the 48 lines have no points over $\mathbf{F}_{47}$ either.

Note that this method is limited since by the Hasse-Weil theorem, $\mid 4a_q \mid \leq 8\sqrt{q}$, and so for large enough $q$ (in fact $q > 61$), $\#X(\mathbf{F}_q) > 0$.

**(10)** $T = \{0, 1, 3, 4, 6, 7, 9\}, t = 6$. Then $V(T,t)$ is the intersection of 2 hypersurfaces. Their (very long) equations can be computed using Mathematica. Using Magma we see experimentally that if $q$ is $1 \pmod 3$, then $V(T,t)$ has $41q^3 - 184q^2 + 406q - 413$ points over $\mathbf{F}_q$. There are 40 predictable components that are hyperplanes of dimension 3. Apparently there is a 41st component that is rational. Macaulay2 is currently seeking to identify this component. Once done, this should show that:

**Corollary.** If $\{0, 1, 3, 4, 6, 7, 9\} \subseteq S$ and 3 does not divide $n$, then $d \geq 7$.

## 4. Theorems and Conjectures.

We give below some general results regarding the varieties.

**Theorem.** Let $T = \{0, 1, 3, r\}$. Then $V(T, 3)$ consists of $r - 2$ points, given by a polynomial $F_r(x_1, x_2) = 0$ and $x_3 = -x_1 - x_2$. If $r$ is odd, then $x_1, x_2, x_1 + x_2$ are all factors of $F_r$. If $r$ is $0 \pmod 3$, then $(x_1^2 + x_1x_2 + x_2^2)^2$ is a factor of $F_r$. If $r$ is $1 \pmod 3$, then $(x_1^2 + x_1x_2 + x_2^2)$ is a factor of $F_r$. Removing the above factors leave a factor of degree $6k$, which is typically irreducible with Galois group depending only on $k$.

*Proof.* If we set $x_1 = 0$, then $\Delta[0, 1, r]$ simplifies to $x_2x_3^r - x_2^r x_3$. Setting $x_3 = -x_1 - x_2 = -x_2$ yields $(-1)^r x_2^{r+1} + x_2^{r+1}$, which is 0 if $r$ is odd. Divisibility by $(x_1^2 + x_1x_2 + x_2^2)$ is handled likewise by setting $x_1 = \omega x_2$. The case of $(x_1^2 + x_1x_2 + x_2^2)^2$ is handled by differentiating and setting $x_1 = \omega x_2$.

*Remark.* $F_{11}$ is $x_1x_2(x_1+x_2)(x_1^2+x_1x_2+2x_2^2)(2x_1^2+x_1x_2+x_2^2)(2x_1^2+3x_1x_2+2x_2^2)$. Otherwise, for $k = 1$, the sextic has Galois group $S_3$ of order 6. For $k = 2$, the Galois group is the group of order 72 described earlier, namely TransitiveGroup(12, 35) in the Magma database. This is in fact $S_3 \wr C_2$. Exploiting the $S_3$ symmetry of $V(T, 3)$, we obtain the following result.

**Theorem.** The Galois group of $F_r(x, 1)$ embeds naturally in $S_3 \wr S_k$ ($S_3$ acting regularly on 6 letters), where $k$ is $[(r-2)/6]$ unless $r$ is 3) mod 6) when it is $(r-9)/6$.

Note that these Galois groups are typically not solvable and the embeddings are typically surjective. Thus there will in general be no succinctly stated corollaries (since there is no simple description of the residue fields). We next turn to the predictable hyperplanes.

**Theorem.** Suppose that $T = \{r, r+1, ..., r+t-1, r+m, r+m+1, ..., r+m+t-1, ..., r+km, ..., r+km+t-1\}$. Then $V(T, t+k)$ contains $(\frac{t+k}{k-1})(m-1)(m-2)...(m-t)$ $k-1$-dimensional linear subspaces.

*Proof.* It is clear that $\Delta[T]$ vanishes if we take a point on any of the following linear spaces. Let $k-1$ of the variables be free (these can be chosen in $(\frac{t+k}{k-1})$ ways). The remaining $t+1$ variables will be set to be $1, \zeta^{i_1}, ..., \zeta^{i_t}$, where $\zeta$ is a primitive $m$th root of 1 and $0, i_1, ..., i_t$ are distinct modulo $m$. Then $i_1, ..., i_t$ can be chosen in $(m-1)...(m-t)$ ways.

**Theorem.** $V(\{0, 1, m, m+1, 2m\}, 4) = V(\{0, 1, m, m+1, 2m, 2m+1\}, 4)$ and in addition to the predicted lines from the previous theorem, it contains a further component of the form $V(\{0, 1, 2, m\}, 4) \cap (\mathbf{P}^1 \times V(\{0, 1, m\}, 3))$.

*Proof.* $V(\{0, 1, m, m+1, 2m, 2m+1\}, 4) \subseteq V(\{0, 1, m, m+1, 2m\}, 4)$. Equality occurs since if $\underline{x} \in V(\{0, 1, m, m+1, 2m\}, 4)$, then each $x_i^{m+1}$ and $x_i^{2m}$ is a linear combination of $1, x_i, x_i^m$. Say $x_i^{m+1} = A + Bx_i + Cx_i^m$ and $x_i^{2m} = D + Ex_i + Fx_i^m$ (note that $A, B, C, D, E, F$ are independent of $i$). Thus $x_i^{2m+1} = Ax_i^m + Bx_i^{m+1} + Cx_i^{2m} = AB + CD + (B^2 + CE)x_i + (A + BC + CF)x_i^m$.

The second part follows similarly. Namely, let $\underline{x} \in V(\{0, 1, 2, m\}, 4) \cap (\mathbf{P}^1 \times V(\{0, 1, m\}, 3))$. For $i \geq 2$, $x_i^m$ is a linear combination (independent of $i$) of $1, x_i$. Also, $x_i^2$ is a linear combination of $1, x_i, x_i^m$. So $x_i^{m+1}$ is a linear combination of $x_i, x_i^2$, so of $1, x_i, x_i^m$. Similarly, $x_i^{2m}$ is a linear combination of $x_i^m, x_i^{m+1}$, which by what we just showed is a linear combination of $1, x_i, x_i^m$.

The only problem now is in dealing with $x_1$, on which being in $\cap(\mathbf{P}^1 \times V(\{0, 1, m\}, 3)$ imposes no conditions. The point here is that if $\Delta[\{0, 1, 2, m\}] = \Delta[\{0, 1, m\}] = 0$, then the matrix with rows $1, 1, 1, x_2, x_3, x_4, x_2^m, x_3^m, x_4^m$ has zero determinant too. This is an identity easily checked.

**Conjecture.** $V(\{0, 1, m, m+1, 2m\}, 4)$ exactly consists of the predicted lines together with the component of the form $V(\{0, 1, 2, m\}, 4) \cap (\mathbf{P}^1 \times V(\{0, 1, m\}, 3))$.

Note that the above theorem implies that for $m \geq 5$, we will run into a component that is not rational, since $V(\{0, 1, m\}, 3)$ is a curve of genus $\geq 1$. This in turn implies that we will not have the succinctly stated results we desire, in this case. A more promising (but perhaps over-ambitious) case is:

**Question.** If $T = \{0, 1, 3, 4, 6, 7, 9, 10, ..., r\}$, then $V(T, \#T - 1)$ is rational over the field $\mathbf{Q}(\omega)$, and if $T \subseteq S$ and 3 does not divide $n$, then $d \geq \#T$.

The cases of $r \leq 7$ have been established in above examples. The case $r = 9$ is being investigated by computer (see example 10 above).

**Conclusions.**

The conclusion is that for *any* set $T$ of integers and positive integer $t$ we get *some* result of the form that if $T \subseteq S$ and some other condition holds, then $d > t$. This other condition is most succinct (e.g. that 3 not divide $n$) when the variety $V(T,t)$ is rational over a cyclotomic field. If rationality fails or the field of definition is more complicated than a cyclotomic field, then we get some result - thus, our method responds to the question of what happens if the conditions are not right for use of BCH, Hartmann-Tzeng,... (as is usually the case with e.g. quadratic residue codes) although the answers are not simply stated.